\documentclass[14pt, twoside, a4paper]{article}
\let\getprepared\relax
\ifx\TestIngCommand\undefined\relax
\else\input ../../../proc/ppreamb-book.tex 
  \input init.tex \fi
\let\TestIngCommand\undefined

\usepackage{amsmath,amscd,amssymb}                                         
\usepackage{graphics}                                                     
\newtheorem{lem}{Lemma}                                                    
\newtheorem{prop}{Proposition}                                             

\newskip\ttglue\ttglue=.5em plus.25em minus.15em                           
\def\firstname#1{\def\FIRSTNAME{#1}\ignorespaces}
\def\lastname#1{\def\LASTNAME{#1}\ignorespaces}
\def\middleinitial#1{\def\MIDDLEINI{#1}\ignorespaces}
\def\department#1{\def\DEPARTMENT{#1}\ignorespaces}
\def\institute#1{\def\INSTITUTE{#1}\ignorespaces}
\def\address#1{\def\ADDRESS{#1}\ignorespaces}
\def\country#1{\def\COUNTRY{#1}\ignorespaces}
\def\otheraffiliation#1{\def\OTHERAFFILIATION{#1}\ignorespaces}
\def\email#1{\def\EMAIL{#1}\ignorespaces}
\newcount\autcount\autcount=0                                              
\newcount\affcount\affcount=0                                              
\newcount\numcount\newcount\nummcount                                      
\newcount\nummmcount\newcount\nummmmcount                                  
\def\writename#1#2{\ \kern-1ex\hbox{
  \csname AUthor\the#1\endcsname\                                          
  \edef\TESTSTR{}\expandafter\ifx\csname auTHor\the#1\endcsname\TESTSTR    
  \else\csname auTHor\the#1\endcsname.\ \fi                                
  \csname authOR\the#1\endcsname$^{\csname AFF\the#1\endcsname}$
  \expandafter\ifx\csname corr\number#1\endcsname\relax                    
  \else\thanks{Corresponding author.}\ \fi                                 
  }\ifnum#1<#2, \else\ \kern-1ex\fi}
\def\writeemail#1{
  \nummcount=0\relax\nummmcount=0\relax                                    
  \loop\ifnum\nummcount<\autcount\advance\nummcount by1\relax              
    {\expandafter\ifnum\csname AFF\the\nummcount\endcsname=#1\relax        
    \global\advance\nummmcount by1\fi}\repeat                              
  \nummcount=0\relax\nummmmcount=0\relax                                   
  \loop\ifnum\nummcount<\autcount\advance\nummcount by1\relax              
    {\expandafter\ifnum\csname AFF\the\nummcount\endcsname=#1\relax        
    \global\advance\nummmmcount by1\relax\def\blank{}\expandafter          
    \ifx\csname EMAIL\the\nummcount\endcsname\blank(no e-mail)
    \else\csname EMAIL\the\nummcount\endcsname                             
    \fi                                                                    
    \ifnum\nummmmcount<\nummmcount; \fi\fi}\repeat}
\long\def\BeginAuthorList#1\EndAuthorList{#1\relax                         
  \author{\vbox{\hsize=390pt\noindent\numcount=0\relax                     
    \loop\ifnum\numcount<\autcount\advance\numcount by1\relax              
      \writename{\numcount}{\autcount}
      \repeat}\\[2mm]                                                      
    \vbox{\small\numcount=0\relax                                          
      \loop\ifnum\numcount<\affcount\advance\numcount by1\relax            
        \vbox{{\count0=\numcount\relax                                     
          \loop\expandafter\ifnum\csname AFF\the\count0\endcsname
            <\numcount\relax\advance\count0 by1\relax\repeat               
          $^{\csname AFF\the\count0\endcsname}$}
        \def\BLANK{}\expandafter\ifx\csname DEPT\the\numcount\endcsname    
          \BLANK                                                           
          \else\csname DEPT\the\numcount\endcsname, \fi                    
        \csname INST\the\numcount\endcsname,                               
        \csname ADDR\the\numcount\endcsname,                               
        \csname COUN\the\numcount\endcsname                                
        \edef\TEST{}\expandafter\ifx\csname OTHE\the\numcount\endcsname
          \TEST                                                            
          .\else;\break\csname OTHE\the\numcount\endcsname.\fi}
        \vbox{\writeemail{\numcount}}
        \repeat}\\}}
\expandafter\def\csname x1\endcsname{}
\expandafter\def\csname x2\endcsname{}
\expandafter\def\csname x3\endcsname{}
\expandafter\def\csname x4\endcsname{}
\expandafter\def\csname x5\endcsname{}
\expandafter\def\csname x6\endcsname{}
\expandafter\def\csname x7\endcsname{}
\expandafter\def\csname x8\endcsname{}
\expandafter\def\csname x9\endcsname{}
\def\Author#1#2{\global\advance\autcount by1\relax#2                       
  \expandafter\edef\csname AUthor\the\autcount\endcsname{\FIRSTNAME}
  \expandafter\edef\csname auTHor\the\autcount\endcsname{\MIDDLEINI}
  \expandafter\edef\csname authOR\the\autcount\endcsname{\LASTNAME}
  \expandafter\edef\csname EMAIL\the\autcount\endcsname{\EMAIL}
  \let\tempera\"\def\"{\string\"}\expandafter\ifx\csname x\DEPARTMENT
    \endcsname\relax                                                       
    \global\advance\affcount by1\relax\let\"\tempera                       
    \expandafter\edef\csname DEPT\the\affcount\endcsname{\DEPARTMENT}
    \expandafter\edef\csname INST\the\affcount\endcsname{\INSTITUTE}
    \expandafter\edef\csname ADDR\the\affcount\endcsname{\ADDRESS}
    \expandafter\edef\csname COUN\the\affcount\endcsname{\COUNTRY}
    \expandafter\edef\csname OTHE\the\affcount\endcsname{\OTHERAFFILIATION}
    \expandafter\edef\csname AFF\the\autcount\endcsname{\the\affcount}
  \else\expandafter\edef\csname AFF\the\autcount\endcsname{\DEPARTMENT}
  \fi\let\"\tempera\ignorespaces}
\def\CorrespondingAuthor#1#2{
  \expandafter\xdef\csname corr\number#1\endcsname{cor}
  \Author#1{#2}}
\def\PaperTitle#1{\title{\bf#1}}
\def\Category#1{\ignorespaces}
\def\keywords#1{{\noindent \emph{Keywords:}                                
  \def\BLANK{}\def\TEST{#1}\ifx\BLANK\TEST(n/a).\else#1\fi}}
\getprepared                                                               
\begin{document}                                                           

\PaperTitle{An Explanation of Mellin's 1921 Paper}
\Category{(Pure) Mathematics}

\date{}

\BeginAuthorList
  \Author1{
    \firstname{Wayne}
    \lastname{Lawton}
    \middleinitial{M}   
    \department{Department of the Theory of Functions, Institute of Mathematics and Computer Science}
    \institute{Siberian Federal University}
    \otheraffiliation{}
    \address{Krasnoyarsk}
    \country{Russian Federation}
    \email{wlawton50@gmail.com}}
\EndAuthorList
\maketitle
\thispagestyle{empty}
\begin{abstract}
In 1921 Mellin published a Comptes Rendu paper computing the principal root of the polynomial $Z^n + x_1Z^{n_1} + \cdots + x_pZ^{n_p} - 1$ using hypergeometric functions of its coefficients $x_1,...,x_p.$  He used an integral transform nowadays bearing his name. Slightly over three pages, the paper is written in French in a terse style befitting the language. Unable to find an elementary explanation on the web or in a texbook, we wrote this expository article to make Mellin's landmark result accessible to interested people who are not experts in hypergeometric functions and complex analysis. 
\end{abstract}
%
%
\footnote{\thanks{This work is supported by the Krasnoyarsk Mathematical 
Center and financed by the Ministry of Science and Higher Education 
of the Russian Federation in the framework of the establishment and 
development of Regional Centers for Mathematics Research and 
Education (Agreement No. 075-02-2020-1534/1).}}
\section{Principal Solution $Z$ and Mellin Transform of $Z^\alpha$}
In the opening paragraph of his paper, Mellin says that it summarizes research he undertook years ago and was prompted by notes  \cite{birkeland} of Richard Birkeland, a Norwegian mathematician known for his contributions to the theory of algebraic equations.
The 7 numbered equations in Mellin's paper and this paper coincide. We prove each of them.
\begin{equation}\label{eqn1}
    Z^n + x_1\, Z^{n_1} + \cdots + x_p\, Z^{n_p} -1 = 0
\end{equation} 
For integers $p \geq 1,$ $0 < n_p < \cdots n_1 < n,$ we define a principal solution to be an analytic function $Z(x_1,...,x_p)$ on $[0,\infty)^p$ satisfying (1) and $Z(0,...,0) = 1.$ 
%
%
\begin{lem}\label{lem1} 
	If a principal solution $Z(x_1,...,x_p)$ of (1) exists, it is unique.
\end{lem}
{\bf Proof}
For $r > 0$ define $D_r := \{z \in \mathbb C:|z|<r\}.$
For sufficiently small $r$ a principal solution $Z(x_1,\dots,x_p)$ of (1) extends to give a holomorphic 
solution $\widetilde Z(x_1,\dots,x_p)$ of (\ref{eqn1}) for $(x_1,...,x_n)$ in the polydisc $D_r^p.$ 
Therefore, for every $n$--th root of unity $\epsilon,$ the function 
$Z_{\epsilon}(x_1,...,x_n) := \epsilon\, \widetilde Z(\epsilon^{n_1}\, x_1,...,\epsilon^{n_p}x_p)$ is a holomorphic solution of (\ref{eqn1}) in $D_r^p.$ Since a polynomial of degree $n$ can have at most $n$ distinct roots, these $n$ distinct functions describe all holomorphic solutions of 
(\ref{eqn1}) on $D_r^p.$ The conclusion follows since the restrictions of $Z(x,...,x_p)$ and $Z_1(x_1,...,x_p)$ to $[0,r)^p$ are equal and analytic.
%
%
Define $\Xi := \{\, \xi := (\xi_1,...,\xi_p) \, : \, \xi_1+\cdots + \xi_p \in (-\infty,-1]\, \},$ $W(\xi) := 1+\xi_1+\cdots+\xi_p,$ and the holomorphic function 
$\Psi := (x_1,...,x_n) : \mathbb C^p \backslash \Xi \rightarrow \mathbb C^p$ by
\begin{equation}\label{eqn2}
    x_i(\xi) := \xi_i\, W(\xi_1,...,\xi_p)^{\frac{n_i}{n} - 1},\, i = 1,...,p
\end{equation}
where $W^\frac{1}{n} : \mathbb C^p \backslash \Xi \rightarrow \mathbb C$ is the branch of $\sqrt[n]{W}$ that satisfies $W^\frac{1}{n}(0,...,0) = 1.$
%
%
\begin{lem}\label{lem2}
If $x_1,...,x_p$ are defined by (\ref{eqn2}), then $Z := W^{-\frac{1}{n}}$ satisfies (\ref{eqn1}). 
\end{lem}
{\bf Proof} Follows by substitution.
%
%
%
\begin{lem}\label{lem3}
The restriction 
$\Psi : [0,\infty)^p \rightarrow  [0,\infty)^p$ is a bijection.
Furthermore, the Jacobian of
$\Psi : \mathbb C^p \backslash \Xi \rightarrow \mathbb C^p$ satisfies
\begin{equation}\label{eqn3}
\frac{\partial(x_1,...,x_p)}{\partial(\xi_1,...,\xi_p)} = 
\left(1+\sum_{k=1}^p \xi_k \right)^{\frac{n_1+\cdots+n_p}{n} - p - 1}
\left(1+\frac{1}{n} \sum_{k=1}^p n_k\, \xi_k \right).
\end{equation}
The restriction $\Psi : [0,\infty)^p \rightarrow  [0,\infty)^p$ and its inverse
$\Psi^{-1} : [0,\infty)^p \rightarrow  [0,\infty)^p$ are analytic.
The principal solution of (\ref{eqn1}) is $Z := W^{-\frac{1}{n}} \circ \Psi^{-1} : [0,\infty)^p \rightarrow \mathbb [1,\infty).$ 
\end{lem}
{\bf Proof}
The first assertion is Proposition \ref{prop1}.
Define the Kronecker symbol
$$
\delta_{i,j} = \begin{cases}
			1, \ \hbox{ if } i = j \\
			0 \ \hbox{ if } i \neq j.
			\end{cases}
$$
Since the matrix
$$\frac{\partial \Psi_i}{\partial \xi_j} = 
W^{\frac{n_i}{n}-1} \left( \delta_{i,j} + \xi_i \, (\frac{n_i}{n}-1)\, W^{-1}\right),$$ 
is a product of two matrices, it follows that
$$\frac{\partial(x_1,...,x_p)}{\partial(\xi_1,...,\xi_p)} := 
\det \frac{\partial \Psi_i}{\partial \xi_j} =  W^{\frac{\sum_{k=1}^p n_k}{n} - p} \, 
\det \left( \delta_{i,j} + \xi_i \, (\frac{n_i}{n}-1)\, W^{-1}\right).$$
Proposition \ref{prop2} implies that
$$\det \left( \delta_{i,j} + \xi_i \, W^{-1}\right) = 
1 + \sum_{k=1}^p \xi_k \, \left(\frac{n_k}{n}-1\right)\, W^{-1} = 
W^{-1}\left(1+\frac{1}{n}\sum_{k=1}^p n_k\,\xi_k \right)
$$
and concludes the proof of the second assertion. Since the Jacobian of $\Psi$ is nozero and holomorphic on $\mathbb C^p \backslash \Xi$ and $[0,\infty)^p \subset \mathbb C^p \backslash \Xi,$ assertion three follows from 
the inverse function theorem for holomorphic functions (\cite{shabat}, p. 40, Theorem 2). The fourth assertion then follows from Lemma \ref{lem2}.
%
%
\\ \\
{\bf Remark 1}
Let $n = 2, p = 1.$ Then 
$\Psi : \mathbb C \backslash (-\infty,-1] \rightarrow 
\mathbb C \backslash i((-\infty,-2] \cup [2,\infty))$
is a holomorphic bijection with inverse
$$
\Psi^{-1}(z) = -1+\left(\, \frac{z}{2}+\sqrt {1+\left(\frac{z}{2}\right)^2}\ \right)^2.
$$
Letting 
$\xi_1 = -1+e^{s+it}$ with 
$s \in \mathbb R$
and
$-\pi < t < \pi$ 
gives
$\Psi(\xi_1) = 2 \sinh \frac{s+it}{2} = u+iv$
where
$u = 2\sinh \frac{s}{2}\cos t/2$
and
$v = 2\cosh \frac{s}{2} \sin t/2.$
For $t = 0$ this gives the curve $v = 0.$ 
For fixed $t \neq 0$ this gives the  
branch of the hyperbola described by the equation
$\left(\frac{v}{ \sin t/2} \right)^2 - 
\left(\frac{u}{\cos t/2} \right)^2 = 4$
with $t v \geq 0.$
%
%
\\ \\
%
%
{\bf Question 1}
For $n \geq 3$ and $p \geq 1,$ is $\Psi$ injective and what is its image?
%
%
\\ \\
Following (\ref{eqn3}) Mellin says that using the known formula, derived in Proposition 2, in combination with 
(\ref{eqn2}) and (\ref{eqn3}), one can deduce the following result:
%
%
\begin{lem}\label{lem4}
The principal solution of (\ref{eqn1}) satisfies 
\begin{equation}\label{eqn4}
\int_0^\infty \cdots \int_0^\infty 
Z(x_1,...,x_p)^\alpha \, 
x_1^{u_1-1}\cdots x_p^{u_p-1}
\, dx_1\cdots dx_p =
\frac{\alpha}{n}\, \frac{\Gamma(u)\Gamma(u_1)\cdots \Gamma(u_p)}{\Gamma(u+u_1+\cdots+u_p+1)},
\end{equation}
where $\alpha > 0,$ and the real parts of $u := \frac{\alpha}{n} -\frac{n_1}{n}u_1 - \cdots -\frac{n_p}{n}u_p,$
and $u_1,...,u_p$ are positive.
\end{lem}
{\bf Proof}
Equation (\ref{eqn2}) and Jensen's inequality implies
$$W^{(u_1n_1+\cdots+u_pn_p)/n} > p^{u_1+\cdots + u_p}\, x_1^{u_1}\cdots x_p^{u_p}.$$
Therefore, since $Z = W^{-\frac{1}{n}}$ and
$\alpha > u_1n_1+\cdots+u_pn_p,$ the integral $I$ in (\ref{eqn4}) exists. 
Proposition \ref{prop3} implies that $F : [0,\infty)^p \rightarrow [0,\infty)^p$ is a bijective diffeomorphism, so $I$ can be expressed as an integral over $[0,\infty)^p$ with respect to the variables $\xi_1,...\xi_p.$ Hence (\ref{eqn2}) and (\ref{eqn3}) give
$$
	I = \int_0^\infty \cdots \int_0^\infty \frac{(1+\frac{n_1}{n}\xi_1 + \cdots + \frac{n_p}{n}\xi_p)
	\xi_1^{u_1-1}\cdots \xi_p^{u_p-1}}{W^\omega}\, d\xi_1\cdots d\xi_p
$$ 
where
$\omega := u + u_1+ \cdots + u_p + 1.$ Proposition (\ref{prop2}) implies that $I + I_0 + I_1+\cdots+I_p$ where
$$I_0 = \frac{\Gamma(\omega-u_1-\cdots - u_p)\Gamma(u_1)\cdots \Gamma(u_p)}{\Gamma(\omega)} = 
\frac{\Gamma(u+1)\Gamma(u_1)\cdots \Gamma(u_p)}{\Gamma(\omega)}
,$$
and
$$I_i = \frac{n_i}{n}\, \frac{\Gamma(\omega-u_1-\cdots - u_p-1)\Gamma(u_1)\cdots \Gamma(u_i+1) \cdots \Gamma(u_p)}{\Gamma(\omega)} = \frac{n_iu_i}{nu}\, I_0.$$
We conclude the proof by observing that
$$I = \left(1 + \frac{1}{u} \sum_{i=1}^p \frac{n_iu_i}{n} \right) \, I_0 = 
\frac{\alpha}{nu}I_0 = \frac{\alpha}{n} \frac{\Gamma(u+1)}{u} 
\frac{\Gamma(u_1)\cdots \Gamma(u_p)}{\Gamma(\omega)}.$$
%
%
\section{Computing $Z^\alpha$ from its Mellin Transform} 
Mellin continues: ``The law of reciprociy relating to the integrals of this species, demonstrated by us in a previous work \cite{mellin1}, allows us to invert as follows (\ref{eqn4}):'' 
\\ \\
$Z(x_1,...,x_p)^\alpha =$
\begin{equation}\label{eqn5}
\frac{1}{(2\pi i)^p} 
\int_{a_1-i\infty}^{a_1+i\infty} \cdots \int_{a_p-i\infty}^{a_p+i\infty} 
\frac{\alpha}{n}\, \frac{\Gamma(u)\Gamma(u_1)\cdots \Gamma(u_p)}
{\Gamma(u+u_1+\cdots + u_p +1)}\, 
x_1^{-u_1}\cdots x_p^{-u_p}\,
du_1\cdots du_p,
\end{equation}
$$\alpha - n_1a_1 -\cdots - n_pa_p > 0, \ a_s > 0, \ 
u := \frac{\alpha}{n} - \frac{n_1}{n} u_1- \cdots -\frac{n_p}{n}u_p$$
$$u+u_1+\cdots + u_p = \frac{\alpha}{n} +\frac{n_1^{`}}{n} u_1+ \cdots 
+\frac{n_p^{`}}{n}u_p, \ \ n_s^{`} := n-n_s.$$
``This formula constitutes our solution of (\ref{eqn1}). It supposes that
$$
-\frac{n_s\pi}{2 n} < \arg (x_s) < \frac{n_s\pi}{2 n}
$$
but we can extend our domain of validity by suitably deforming the integration paths.''
\\ \\
The right side of (\ref{eqn4}) is the Mellin transform of $Z(x_1,...,x_p)^{\alpha}$ and (\ref{eqn5}) represents $Z(x_1,...,x_p)^{\alpha}$ by the inverse Mellin transform.
Lacking the luxury of accessing \cite{mellin1} we refer the 
reader to the derivation of the inverse Mellin transform via the
Fourier transform by Debnath and Bhatta in section 8.1 of \cite{debnathbhatta}. They address the univariate case, but 
the extension to multivariate case is straightforward.
\section{Generalized Hypergeometric Functions}
Let $n$ be a positive integer and $f_1,...,f_p, g_1,...,g_s : \mathbb C^p \rightarrow \mathbb C$ be entire functions and $F : \mathbb C^p \rightarrow \mathbb C$ satisfy the following system of  functional equations
\begin{equation}\label{eqn6}
F(u_1,...,u_s+n,...,u_p) = \frac{f_s(u_1,...,u_p)}{g_s(u_1,...,u_p)}\, F(u_1,...,u_p), \ \ s = 1,...,p.
\end{equation}
and such that the following integral converges and does not change when we move the (vertical) integration paths for each $u_s$ $n$ units to the right:
$$
y(x_1,...,x_p) = \frac{1}{(2\pi i)^p}\, 
\int_{u_1} \cdots \int_{u_p}, F(u_1,...,u_p)\, x_1^{-u_1}\cdots x_p^{-u_p}, \ \ (x_1,...,x_p) \in [0,\infty)^p.
$$
Mellon cites a result \cite{mellin1}
that $y$ satisfies the following system of partial differential equations
\begin{equation}\label{eqn7}
f_s\left(-x_1\frac{\partial}{\partial x_1},...,
-x_p\frac{\partial}{\partial x_p}\right)\, y
= g_s\left(-x_1\frac{\partial}{\partial x_1},...,
-x_p\frac{\partial}{\partial x_p}\right)\, s_s^n\, y, \ \ s = 1,...,p.
\end{equation}
{\bf Remark 2} (\ref{eqn7}) is a system of PDE's of finite order iff $f_s$ and $g_s$ are polynomials.
\\ \\
Mellin calls solutions of (\ref{eqn6}) hypergeometric type if the factors of $f_s, g_s$ have the form
$$
	c_1u_1+ \cdots + c_pu_p+a
$$
where each $c_i$ is a rational real number.
\\ \\
{\bf Remark 3} Polynomials are not hypergeometric functions or series as defined in \cite{kauers,slater} but they are of hypergeometric type since $\prod_{k=0}^{m-1} (x\frac{d}{dx} - k) = x^m\frac{d^m}{dx^m}.$
\\ \\
Formally (\ref{eqn7}) follows from (\ref{eqn6}) since the functions $x_s^{-u_s}$ are eigenfunctions with eigenvalue $-u_s$ of the Euler operator $x_s\frac{\partial}{\partial x_s}, s  = 1,...,p.$ Clearly the function
$$
F(u_1,...,u_p) = \frac{\alpha}{n}\, \frac{\Gamma(u)\Gamma(u_1)\cdots \Gamma(u_p)}
{\Gamma(u+u_1+\cdots + u_p +1)}
$$
satisfies (\ref{eqn6}) where $f_s$ and $g_s$ have the form above, hence $Z(u_1,...,u_p)^\alpha$ is of hypergeometric type whenever
$\alpha > n_1n.$ 
\\ \\
%
%
{\bf Question 2} How large must $\alpha$ be to ensure that $Z^\alpha$ is a solution of (\ref{eqn1}) of hypergeometric type? The condition $\alpha > nn_1$ is sufficient but not necessary because for $n = 1, p = 1$ the root
$Z(x_1) := -\frac{x_1}{2} + \sqrt {1+(x_1/2)^2}$ is of hypergeometric type. Semusheva and Tsikh \cite{semushevatsikh} proved this fact directly by deriving the following Mellin--Barnes integral representation
$$
	Z(x_1) = \frac{1}{4\pi i} \int_{\Re z = \frac{1}{2}} 
	\frac{\Gamma(z)\, \Gamma((1+z)/2)}{\Gamma((3+z)/2)}\, x_1^{z}\, dz.
$$
{\bf Remark 4} For a comprehensive development of Mellin's solution of (\ref{eqn3}) and
systems of differential equations of hypergeometric type see (\cite{belardinelli}, Chapitre V).
\section{Appendix: Crucial Propositions}
We prove results required to derive the equations in Mellin's paper.
\\ \\
Let $e_1 := (1,0,...,0), e_2 := (0,1,...,0),...,e_p = (0,...,0,1)$ be the standard basis for $\mathbb R^p.$ For $s > 0$ let $H_s$ denote the convex hull of $\{se_i:i=1,...,p\}$ and let
$K_s$ denote the convex hull of $\{s(1+s)^{\frac{n_i}{n}-1}e_i:i = 1,...,p\}.$
%
%
\begin{prop}\label{prop1}
$\Psi : [0,\infty)^p \rightarrow [0,\infty)^p$ is a bijection.
\end{prop}
%
%
{\bf Proof}
For every $s > 0$ define the linear map $L_s : \mathbb R^p \rightarrow \mathbb R^p$ by
$$(L_s \, y)_i : = y_i\, (1+s)^{\frac{n_i}{n}-1}, \ \ i = 1,...,p.$$
If $(\xi_1,...,\xi_p) \in H_s,$ then $W(\xi_1,...,\xi_p) = 1+s$ so Equation \ref{eqn2} implies that
$$\Psi(\xi_1,...,\xi_p) = L_s (\xi_1,...,\xi_p) = \sum_{i=1}^p \frac{\xi_i}{s}\, s(1+s)^{\frac{n_i}{n}-1}e_i \in K_s,$$
hence the restriction $\Psi : H_s \rightarrow K_s$ is a linear bijection.
Clearly $\Psi(0,...,0) = (0,...,0)$ 
and $[0,\infty)^p$ is a disjoint union of $\{(0,...,0)\}$ and the sets $H_s, s > 0.$ Since for $i = 1,...,p,$ the function $s(1+s)^{\frac{n_i}{n}-1}$ is increasing, it follows that
$[0,\infty)^p$ is a disjoint union of $\{(0,...,0)\}$ and the sets $K_s, s > 0.$ This concludes the proof.
%
%
%
%
\begin{prop}\label{prop2}
Let $p \geq 1,$ $y_1,...,y_p$ be indeterminates, and define 
the matrix
$$M(y_1,...,y_p) := \left[
  \begin{array}{ccccc}
    1+y_1 & y_1 & \cdots & y_1 & y_1 \\
    y_2 & 1+y_2 & \cdots & y_2 & y_2 \\
    \vdots & \vdots & \ddots& \vdots & \vdots \\
    y_{p-1} & y_{p-1} & \cdots & 1+y_{p+1} & y_{p+1} \\
    y_p & y_p & \cdots & y_p & 1+y_p \\
  \end{array}
\right].$$
Then
$$\det M(y_1,...,y_p) = 1+y_1+\cdots+y_p.$$
\end{prop}
{\bf Proof}
For $i = 1,...,p$ let $\partial_i := \frac{\partial}{\partial y_i}.$
Clearly $\det M(y_1,...,y_p)$ is a polynomial with constant term 
$\det M(0,...,0) = \det I_p = 1$ and of degree at most $1$ in each variable $y_1,...,y_p$ so for $i = 1,...,p$
$$
\partial_i^2 \det M(y_1,...,y_p) = 0.
$$
Let $M_i(y_1,...,y_p)$ denote the matrix obtained from $M(y_1,...,y_p)$ by replacing each entry in its $i$--th row with $1$ and for $j \neq i$ let
$M_{i,j}(y_1,...,y_p)$ denote the matrix obtained from $M(y_1,...,y_p)$
by replacing each entry in its $i$--row and its $j$--row by $1.$
Since a determinant of a matrix is a linear function of each of its row vectors, 
$$\partial_i \det M(y_1,...,y_p) := \det M_i(y_1,...,y_p)$$
and for $j \neq i$
$$\partial_j \partial_i \det M(y_1,...,y_p) := \det M_{i,j}(y_1,...,y_p) = 0$$
since $M_{i,j}(y_1,...,y_p)$ has two identical rows.
Therefore Taylor's expansion gives
$$\det M(y_1,...,y_n) = 1 + \sum_{i = 1}^p y_i \, \partial_i \det M(0,...,0)
= 1 + \sum_{i = 1}^p y_i \, \det M_i(0,...,0).$$
It suffices to prove
$\det M_i(0,...,0) = 1, \, i = 1,...,p.$
This follows since
$$ M_i(0,...0) =  \left[
\begin{array}{ccccccc}
    1 & 0 & \cdots & 0 & \cdots & 0 & 0 \\
    0 & \ddots & \ddots & \ddots & \ddots & \ddots & 0\\
    0 & \cdots & 1 & 0 & 0 & \cdots & 0\\
    1 & \cdots & 1 & 1 & 1 & \cdots & 1\\
    0 & \cdots & 0 & 0 & 1 & \cdots & 0 \\
    0 & \ddots & \ddots & \ddots & \ddots & \ddots & 0\\
    0 & 0 & \cdots & 0 & \cdots & 0 & 1 \\
  \end{array}
\right].$$
%
%
{\bf Remark 5}
Nilsson and Passare (\cite{nilssonpassare}, Example 1) used basic calculus to compute the following integrals for 
$\Re u_i \in (0,1), i = 1,2.$
$$
	\int_0^\infty \int_0^\infty 
	\frac{\xi_1^{u_1-1}\, \xi_2^{u_2-1}}{1+\xi_1+\xi_2} d\xi_1 d\xi_2 = 
	\frac{\Gamma(u_1)\, \Gamma(u_2) \, \Gamma(1-u_1-u_2)}{\Gamma(1)}
$$
The following result uses exterior calculus to extend their computation.
%
%
\begin{prop}\label{prop3}
For every positive integer $p,$ complex $u_1,...,u_p$ satisfying $\Re \, u_i > 0,$ and $\omega > \max \{\Re \, u_i\}$ the integral below converges and satisfies the stated identity.
$$
\int_0^\infty \cdots \int_0^\infty 
\frac{\xi_1^{u_1-1}\cdots \xi_p^{u_p-1}}{\left( 1+\xi_1+\cdots+\xi_p \right)^\omega}\, d\xi_1\cdots d\xi_p =
\frac{\Gamma(u_1)\cdots \Gamma(u_p)\Gamma(\omega-u_1-\cdots - u_p)}{\Gamma(\omega)}.
$$
\end{prop}
{\bf Proof}
Multiplying both sides of the asserted identity in the lemma by 
$$\Gamma(\omega) = \int_0^\infty x^{\omega-1}\, e^{-x}\, dx$$
gives the equivalent identity 
$
I_1(u_1,...,u_p,\omega) = 
\Gamma(u_1)\cdots \Gamma(u_p)\, \Gamma(\omega-u_1-\cdots - u_p)
$ 
where
$$
I_1(u_1,...,u_p,\omega) := 
\int_0^\infty \cdots \int_0^\infty 
\frac{\xi_1^{u_1-1}\cdots \xi_p^{u_p-1}}{\left( 1+\xi_1+\cdots+\xi_p \right)^\omega}\, x^{\omega - 1}\, e^{-x}\, d\xi_1\wedge \cdots \wedge d\xi_p \wedge dx.
$$
Here we express the volume form using exterior products to
implement the following change of variables: 
$y := x/(1+\xi_1+\cdots+\xi_p)$ and $z_i := y\, \xi_i, i = 1,...,p.$ 
Then $I_1(\xi_1,...,\xi_p,x) = I_2(z_1,...,z_2,y)$ where
$$
I_2(z_1,...,z_2,y) := \int_0^\infty \cdots \int_0^\infty 
z_1^{u_1-1}e^{-z_1}\cdots z_p^{u_p-1}e^{-z_p}\, y^{\omega - u_1-\cdots -u_p-1}\, e^{-y}\, \theta
$$
where 
$$
	\theta = \frac{y^p}{1+\xi_1+\cdots+\xi_p}\, d\xi_1\wedge \cdots \wedge d\xi_p \wedge dx
$$
In order to finish the proof it suffices to prove that $dz_1 \wedge \cdots \wedge dz_p \wedge dy = \theta.$
Since $dz_i = \xi_i dy + y\, d\xi_i$ and $dy \wedge dy = 0$ it follows that
$$dz_1 \wedge \cdots \wedge dz_p \wedge dy = y^p\, d\xi_1 \wedge \cdots \wedge d\xi_p \wedge dy.$$
Substitute
$$
	dy = \frac{1}{1+\xi_1+\cdots+\xi_p}\, dx - \frac{x}{(1+\xi_1+\cdots+\xi_p)^2} \left(\sum_{k=1}^p d\xi_k\right)
$$
and use the fact that $d\xi_i \wedge d\xi_i = 0$ to obtain
$$
	dz_1 \wedge \cdots \wedge dz_p \wedge dy = \frac{y^p}{1+\xi_1+\cdots+\xi_p}\, d\xi_1\wedge \cdots \wedge d\xi_p \wedge dx = \theta 
$$
which concludes the proof.
%
%
\\ \\
{\bf Acknowledgments} We thank Elaine Wong of the Johann Radon Institute for Computational and Applied Mathematics (RICAM), Linz for her meticulous proofreading of and corrections to this article.


\begin{thebibliography}{10}

\bibitem{belardinelli} G. Belardinelli, {\it Fonctions hyperg\'eom\'etriques de plusieurs variables et r\'esolution analytique des \'equations alg\'ebriques g\'en\'erales.} M\'emorial des sciences math\'ematiques, no. 145, Gauthier--Villars, Paris, 1960.

\bibitem{birkeland} M. R. Birkeland, {\it Ein allgemeiner Satz \"uber algebraische Gleichungen.} Annales Academiae Scientarum Fennicae, ser. A, t. 7, (1915).

\bibitem{debnathbhatta} L. Debnath and D. Bhatta,
{\it Integral Transforms and Their Applications, 2nd edition.},
Chapman \& Hall/CRC, New York, 2007.

\bibitem{kauers} M. Kauers and P. Paule, {\it The Concrete Tetrahedron, Symbolic Sums, Recurrence Equations, Generating Functions, Asymptotic Estimates.} Springer, Vienna, 2011.

\bibitem{mellin1} H. J. Mellin, 
{\it Zur Theorie zweier allgemeiner Klassen bestimmter Integrale.} 
75 S. Ibid, Bd. 22 (1896).

\bibitem{mellin2} H. J. Mellin, {\it R\'esolution de l'equation 
alg\'ebrique g\'en\'erale \'a  
l'aide de la fonction gamma.} 
C. R. Acad. Sci. Paris S\'er. I Math. 172 (1921) 658--661.
%
\bibitem{nilssonpassare} L. Nilsson and M. Passare, 
{\it Mellin transforms of multivariate rational functions.} 
Journal of Geometric Analysis, 23 (2013) 24--46.
arXiv:1010.5060

\bibitem{semushevatsikh} A. Y. Semusheva and A. K. Tsikh, {\it Continuation
on Mellin's research on solving algebraic equations.} Complex analysis and differential operators. To the 150th anniversary of S. V. Kovalevskaya. Executive editor A. K. Tsikh, Krasnoyarsk State University, Krasnoyarsk, 2000. 196 p. (Russian)

\bibitem{shabat} B. V. Shabat, 
{\it Introduction to Complex Analysis, Part II Functions of Several Variables.}
American Mathematical Society, Providence, Rhode Island, 1992.

\bibitem{slater} L. J. Slater, {\it Generalized Hypergeometric Functions.} 
Cambridge University Press, 1960.

\end{thebibliography}
\end{document}